\documentclass[10pt,aps,prb,reqno]{amsart}
\usepackage{amsmath}
\usepackage{amssymb}

\newtheorem{theorem}{Theorem}[section]
\newtheorem{remark}{Remark}[section]
\newtheorem{lemma}[theorem]{Lemma}
\newtheorem{proposition}[theorem]{Proposition}

\begin{document}
\title[Generalized MHD system]{On the global regularity of two-dimensional generalized magnetohydrodynamics system}
\author{Kazuo Yamazaki}  
\date{}
\thanks{The author expresses gratitude to Professor Jiahong Wu and Professor David Ullrich for their teaching.
}

\maketitle 

\begin{abstract}
We study the two-dimensional generalized magnetohydrodynamics system with dissipation and diffusion in terms of fractional Laplacians. In particular, we show that in case the diffusion term has the power $\beta = 1$, in contrast to the previous result of $\alpha \geq \frac{1}{2}$, we show that $\alpha > \frac{1}{3}$ suffices in order for the solution pair of velocity and magnetic fields to remain smooth for all time.
\vspace{5mm}

\end{abstract}
\footnote{2000MSC : 35B65, 35Q35, 35Q86}
\footnote{Department of Mathematics, Oklahoma State University, 401 Mathematical Sciences, Stillwater, OK 74078, USA}

\section{Introduction and statement of results}

We study the generalized magnetohydrodynamics (MHD) system defined as follows: 

\begin{equation}
\begin{cases}
\frac{\partial u}{\partial t} + (u\cdot\nabla) u - (b\cdot\nabla)b + \nabla \pi + \nu \Lambda^{2\alpha}u = 0,\\
\frac{\partial b}{\partial t} + (u\cdot\nabla) b - (b\cdot\nabla)u + \eta \Lambda^{2\beta} b = 0,\\
\nabla \cdot u = \nabla\cdot b = 0, \hspace{5mm} u(x,0) = u_{0}(x), b(x,0) = b_{0}(x),
\end{cases}
\end{equation}
where $u: \mathbb{R}^{N}\times \mathbb{R}^{+}\mapsto \mathbb{R}^{N}$ is the velocity vector field, $b: \mathbb{R}^{N}\times \mathbb{R}^{+}\mapsto \mathbb{R}^{N}$ the magnetic vector field, $\pi: \mathbb{R}^{N}\times \mathbb{R}^{+}\mapsto \mathbb{R}$ the pressure scalar field and $\nu, \eta \geq 0$ are the kinematic viscosity and diffusivity constants respectively. We also denote by $\Lambda$ a fractional Laplacian operator defined via Fourier transform as $\widehat{\Lambda^{2\gamma}f}(\xi) = \lvert \xi\rvert^{2\gamma} \hat{f}(\xi)$ for any $\gamma \in \mathbb{R}$. 

In case $N = 2, 3, \nu, \eta > 0, \alpha = \beta = 1$, it is well-known that (1) possesses at least one global $L^{2}$ weak solution; in case $N =2$, it is also unique (cf. [19]). Moreover, in any dimension $N \geq 2$, the case $\nu, \eta > 0$, the lower bounds on the powers of the fractional Laplacians at $\alpha \geq \frac{1}{2} + \frac{N}{4}, \beta \geq \frac{1}{2} + \frac{N}{4}$ imply the existence of the unique global strong solution pair (cf. [26]). 

Some numerical study have shown that the velocity vector field may play relatively important role in regularizing effect (e.g. [8], [18]). Starting from the works of [9] and [34], we have seen various regularity criteria of the MHD system in terms of only the velocity vector field (e.g. [1], [4], [6], [7], [10], [25], [28], [33]). Moreover, motivated by the work of [20], the author in [24] showed that in case $N \geq 2, \nu, \eta > 0, \alpha \geq \frac{1}{2} + \frac{N}{4}, \beta > 0$ such that $\alpha + \beta \geq 1 + \frac{N}{2}$, the system (1) even in logarithmically super-critical case still admits a unique global strong solution pair. The endpoint case $\nu > 0, \eta = 0, \alpha = 1 + \frac{N}{2}$ was also completed recently in [23] and [29] (cf. [27] for further generalization). 

On the other hand, in case $N =2$, it is well-known that the Euler equations, the Navier-Stokes system with no dissipation, still admits a unique global strong solution. This is due to the conservation of vorticity which follows upon taking a curl on the system. In the case of the MHD system, upon taking a curl and then $L^{2}$-estimates of the resulting system, every non-linear term has $b$ involved. Exploiting this observation and divergence-free conditions, the authors in [2] showed that in case $N = 2$, full Laplacians in both dissipation and magnetic diffusion are not necessary for the solution pair to remain smooth; rather, only a mix of partial dissipation and diffusion in the order of two derivatives suffices. 

Very recently, the authors in [22] have shown that in case $N  =2$, the solution pair remains smooth in any of the following three cases:

\begin{enumerate}
\item $\alpha \geq \frac{1}{2}, \beta \geq 1$,
\item $\alpha \geq 2, \beta = 0$,
\item $\frac{1}{2} > \alpha \geq 0, 2\alpha + \beta > 2$.
\end{enumerate}
In particular, their result implies that in the range of $\alpha \in [0, \frac{1}{2})$, $\beta$ must satisfy 

\begin{equation}
\beta > 2-2\alpha. 
\end{equation}
These results implied that if $\alpha = 0$, then $\beta > 2$ was necessary to obtain global regularity result. This was improved in [31] to show that either of the following conditions suffices:

\begin{enumerate}
\item $\alpha = 0, \beta > \frac{3}{2}$,
\item $\frac{1}{2} > \alpha > 0, \frac{3}{2} \geq \beta > \frac{5}{4}, \alpha + 2\beta > 3$.
\end{enumerate}
In particular, this implies that in the range of $\alpha \in (0, \frac{1}{2})$, $\beta$ must satisfy 

\begin{equation}
 \beta > \frac{3-\alpha}{2} 
 \end{equation} 
(cf. also [32]). In this paper we make further improvement in this direction. Let us present our results. 

\begin{theorem}
Let $N = 2, \nu, \eta > 0, \alpha > \frac{1}{3}, \beta = 1$. Then for all initial data pair $(u_{0}, b_{0}) \in H^{s}(\mathbb{R}^{2}) \times H^{s}(\mathbb{R}^{2}), s \geq 3$, there exists a unique global strong solution pair $(u,b)$ to (1) such that 

\begin{eqnarray*}
&&u \in C([0,\infty); H^{s}(\mathbb{R}^{2}))\cap L^{2}([0,\infty); H^{s+\alpha} (\mathbb{R}^{2})),\\
&&b \in C([0,\infty); H^{s}(\mathbb{R}^{2}))\cap L^{2}([0,\infty); H^{s+1}(\mathbb{R}^{2})).
\end{eqnarray*}

\end{theorem}

\begin{theorem}

Let $N = 2, \nu, \eta > 0,\alpha \in (0, \frac{1}{3}], \beta \in (1, \frac{3}{2}]$ such that 

\begin{equation}
3 < 2\beta + \frac{2\alpha}{1-\alpha}.
\end{equation}
Then for all initial data pair $(u_{0}, b_{0}) \in H^{s}(\mathbb{R}^{2})\times H^{s}(\mathbb{R}^{2}), s \geq 3$, there exists a unique global strong solution pair $(u, b)$ to (1) such that 

\begin{eqnarray*}
&&u \in C([0,\infty); H^{s}(\mathbb{R}^{2})) \cap L^{2}([0,\infty); H^{s+\alpha}(\mathbb{R}^{2})),\\
&&b \in C([0,\infty); H^{s}(\mathbb{R}^{2})) \cap L^{2}([0,\infty) ; H^{s+\beta}(\mathbb{R}^{2})).
\end{eqnarray*}

\end{theorem}

\begin{remark}
\begin{enumerate}
\item We observe that (4) is equivalent to 

\begin{equation*}
\frac{3}{2} - \frac{\alpha}{1-\alpha} < \beta,
\end{equation*}
and this is a better lower bound than that of (2) or (3) for $\alpha \in (0, \frac{1}{3}]$. 

\item Theorem 1.1 also represents the smaller lower bound for the sum of $\alpha + \beta$ at $1+\frac{1}{3}$ required for the solution pair to remain smooth for all time in comparison to the previous works such as [24] and [26] at $\alpha + \beta \geq 1 + \frac{N}{2}$ in $N$-dimension and [22] at $\alpha + \beta \geq \frac{3}{2}$ in two-dimension. 

\item There are various spaces of functions in which one may obtain local well-posedness of the MHD system. We chose to state above for simplicity. The local theory may be obtained by using mollifiers as done in [14] and we omit the details referring interested readers to [2] where the authors considered (1) in case $N = 2, \nu = 0, \eta > 0, \beta = 1$ and showed in particular the existence of its weak solution pair (cf. also [19] and [26]). 

\item After this work was completed, this direction of research has caught much attention from many mathematicians and a remarkable development with new results has been seen. In particular, we mention that in [3] and [11], the authors obtained the global regularity result in the case $\alpha = 0, \beta > 1$. We also mention numerical analysis results obtained in [21] concerning the interesting case $\alpha = 0, \beta = 1$. 

\end{enumerate}

\end{remark}

In the following section, let us set up notations and summarize key lemmas that will be used repeatedly. Thereafter, we prove our theorems. 

\section{Preliminaries}

Let us denote a constant that depends on $a, b$ by $c(a,b)$ and when the constant is not of significance, let us write $A \lesssim B, A \approx B$ to imply that there exists some constant $c$ such that $A \leq cB, A = cB$ respectively. We also denote partial derivatives and vector components as follows:

\begin{equation*}
\frac{\partial}{\partial t} = \partial_{t}, \hspace{5mm} \frac{\partial}{\partial x} = \partial_{1}, \hspace{5mm} \frac{\partial}{\partial y} = \partial_{2}, \hspace{5mm} u = (u_{1}, u_{2}), \hspace{5mm} b = (b_{1}, b_{2}).
\end{equation*}
For simplicity we also set 

\begin{equation}
w = \nabla \times u, \hspace{3mm} j = \nabla \times b,
\end{equation}

\begin{equation*}
X(t) = \lVert w(t)\rVert_{L^{2}}^{2} + \lVert j(t)\rVert_{L^{2}}^{2}, \hspace{3mm} Y(t) = \lVert \nabla w(t)\rVert_{L^{2}}^{2} + \lVert \nabla j(t)\rVert_{L^{2}}^{2}.
\end{equation*}
 
We use the following well-known inequalities: 

\begin{lemma}
Let $f$ be divergence-free vector field such that $\nabla f \in L^{p}, p \in (1,\infty)$. Then the following inequality holds:

\begin{equation*}
\lVert \nabla f\rVert_{L^{p}} \leq c(p)\lVert \text{curl }  f\rVert_{L^{p}}.
\end{equation*}
\end{lemma}

\begin{lemma}
(cf. [13]) Let $f,g$ be smooth such that $\nabla f \in L^{p_{1}}, \Lambda^{s-1}g \in L^{p_{2}}, \Lambda^{s}f \in L^{p_{3}}, g \in L^{p_{4}}, p \in (1,\infty), \frac{1}{p} = \frac{1}{p_{1}}+\frac{1}{p_{2}} = \frac{1}{p_{3}} + \frac{1}{p_{4}}, p_{2}, p_{3} \in (1, \infty), s > 0.$ Then the following inequality holds:

\begin{equation*}
\lVert \Lambda^{s}(fg) - f\Lambda^{s}g\rVert_{L^{p}} \lesssim (\lVert \nabla f\rVert_{L^{p_{1}}}\lVert \Lambda^{s-1}g\rVert_{L^{p_{2}}} + \lVert \Lambda^{s}f\rVert_{L^{p_{3}}}\lVert g\rVert_{L^{p_{4}}}).
\end{equation*}

\end{lemma}

\begin{lemma}
(cf. [5], [12])
For any $\alpha \in [0, 1], x \in \mathbb{R}^{N}, \mathbb{T}^{N}$ and $f, \Lambda^{2\alpha}f \in L^{p}, p \geq 2$, 

\begin{equation*}
2\int \lvert\Lambda^{\alpha}(f^{\frac{p}{2}})\rvert^{2}dx \leq p \int\lvert f\rvert^{p-2}f\Lambda^{2\alpha}fdx.
\end{equation*}

\end{lemma}

Finally, the following product estimate has proven to be useful (e.g. [15], [16], [17], [30]): 

\begin{lemma}
Let $\sigma_{1}, \sigma_{2} < 1, \sigma_{1} + \sigma_{2} > 0$. Then there exists a constant $c(\sigma_{1}, \sigma_{2}) > 0$ such that 

\begin{equation*}
\lVert fg\rVert_{\dot{H}^{\sigma_{1} + \sigma_{2} - 1}} \leq c(\sigma_{1}, \sigma_{2}) \lVert f\rVert_{\dot{H}^{\sigma_{1}}}\lVert g\rVert_{\dot{H}^{\sigma_{2}}},
\end{equation*}
for $f \in \dot{H}^{\sigma_{1}}(\mathbb{R}^{2}), g\in \dot{H}^{\sigma_{2}}(\mathbb{R}^{2})$. 

\end{lemma}

\section{Proof of Theorem 1.1}

Throughout this section, we assume $\alpha \in (\frac{1}{3}, \sqrt{2} - 1)$ as the case $\alpha \in [\sqrt{2} - 1, \frac{1}{2})$  may be done via slight modification using Gagliardo-Nirenberg inequalities. We note that the restriction of this range of $\alpha$ in particular becomes crucial at (9); we chose the statements of Propositions 3.1-3.3 for simplicity of presentation. We work on 

\begin{equation}
\begin{cases}
\partial_{t} u + (u\cdot\nabla) u - (b\cdot\nabla) b + \nabla \pi + \Lambda^{2\alpha} u = 0,\\
\partial_{t}b + (u\cdot\nabla) b - (b\cdot\nabla) u + \Lambda^{2}b = 0.
\end{cases}
\end{equation}
Taking $L^{2}$-inner products of (6) with $u$ and $b$ respectively, we can get 

\begin{equation}
\sup_{t\in [0,T]}\lVert u(t)\rVert_{L^{2}}^{2} + \lVert b(t)\rVert_{L^{2}}^{2} + \int_{0}^{T} \lVert \Lambda^{\alpha} u\rVert_{L^{2}}^{2} + \lVert \Lambda b\rVert_{L^{2}}^{2} d\tau \leq c(u_{0}, b_{0}, T).
\end{equation}

It has been shown that the following proposition can be attained as long as $\beta \geq 1$ (cf. [22], [31]). We sketch its proof for completeness. 

\begin{proposition}
Let $N = 2, \nu, \eta > 0, \alpha \in (\frac{1}{3}, \sqrt{2} - 1), \beta = 1$. Then for any solution pair $(u,b)$ to (1) in $[0,T]$ there exists a constant $c(u_{0}, b_{0}, T)$ such that   

\begin{equation*}
\sup_{t\in [0,T]}\lVert w(t)\rVert_{L^{2}}^{2} + \lVert j(t)\rVert_{L^{2}}^{2} + \int_{0}^{T}\lVert \Lambda^{\alpha} w\rVert_{L^{2}}^{2} + \lVert \Lambda j\rVert_{L^{2}}^{2} d\tau \leq c(u_{0}, b_{0}, T).
\end{equation*}

\end{proposition}

\proof{

Taking curls on (6), we obtain 

\begin{equation}
\begin{cases}
\partial_{t} w + \Lambda^{2\alpha} w = -(u\cdot\nabla) w + (b\cdot\nabla) j,\\
\partial_{t} j + \Lambda^{2} j = -(u\cdot\nabla) j + (b\cdot\nabla) w + 2[\partial_{1}b_{1} (\partial_{1}u_{2} + \partial_{2}u_{1}) - \partial_{1}u_{1} (\partial_{1}b_{2} + \partial_{2}b_{1})].
\end{cases}
\end{equation}
Taking $L^{2}$-inner products with $w$ and $j$ respectively and using incompressibility of $u$ and $b$, we estimate 

\begin{eqnarray*}
&&\frac{1}{2}\partial_{t} (\lVert w\rVert_{L^{2}}^{2} + \lVert j\rVert_{L^{2}}^{2}) + \lVert \Lambda^{\alpha} w\rVert_{L^{2}}^{2} + \lVert \Lambda j\rVert_{L^{2}}^{2}\\
&=& 2\int [\partial_{1}b_{1} ( \partial_{1}u_{2} + \partial_{2}u_{1}) - \partial_{1}u_{1}(\partial_{1}b_{2} + \partial_{2}b_{1})]j\\
&\lesssim& \lVert \nabla b\rVert_{L^{4}}\lVert \nabla u\rVert_{L^{2}}\lVert j\rVert_{L^{4}}\\
&\lesssim& \lVert j\rVert_{L^{2}}\lVert \nabla j\rVert_{L^{2}}\lVert w\rVert_{L^{2}}\\
&\leq& \frac{1}{2} \lVert \Lambda j\rVert_{L^{2}}^{2} + c\lVert j\rVert_{L^{2}}^{2} \lVert w\rVert_{L^{2}}^{2}
\end{eqnarray*}
by H$\ddot{o}$lder's, Gagliardo-Nirenberg and Young's inequalities. Absorbing diffusive term, (7) and Gronwall's inequality complete the proof of Proposition 3.1. 

}

Next two propositions are the keys to the improvement from previous results:

\begin{proposition}
Let $N = 2, \nu, \eta > 0, \alpha \in (\frac{1}{3}, \sqrt{2} - 1), \beta = 1$. Then for any solution pair $(u,b)$ to (1) in $[0,T]$, for any $\gamma \in (1, 1+\alpha)$, there exists a constant $c(u_{0}, b_{0}, T)$ such that   

\begin{equation*}
\sup_{t\in [0,T]}\lVert \Lambda^{\gamma} b(t)\rVert_{L^{2}}^{2} + \int_{0}^{T} \lVert \Lambda^{1+\gamma} b\rVert_{L^{2}}^{2} d\tau \leq c(u_{0}, b_{0}, T).
\end{equation*}

\end{proposition}

\proof{

We fix $\gamma \in (1, 1+\alpha)$. From the magnetic field equation of (6), we estimate after multiplying by $\Lambda^{2\gamma} b$ and integrating in space as follows:

\begin{eqnarray*}
&&\frac{1}{2}\partial_{t} \lVert \Lambda^{\gamma} b\rVert_{L^{2}}^{2} + \lVert \Lambda^{1+\gamma} b\rVert_{L^{2}}^{2}\\
&\leq& (\lVert (u\cdot\nabla) b\rVert_{\dot{H}^{\gamma -1}} \lVert \Lambda^{1+\gamma} b\rVert_{L^{2}} + \lVert (b\cdot\nabla) u\rVert_{\dot{H}^{\gamma -1}} \lVert \Lambda^{1+\gamma} b\rVert_{L^{2}})\\
&\leq& \frac{1}{2} \lVert \Lambda^{1+\gamma} b\rVert_{L^{2}}^{2} + c(\lVert (u\cdot\nabla) b\rVert_{\dot{H}^{\gamma -1}}^{2} + \lVert (b\cdot\nabla) u\rVert_{\dot{H}^{\gamma -1}}^{2})
\end{eqnarray*}
by H$\ddot{o}$lder's and Young's inequalities. Now we use Lemmas 2.4 and 2.1 to estimate 

\begin{eqnarray*}
\lVert (b\cdot\nabla) u\rVert_{\dot{H}^{\gamma -1}}^{2} \lesssim \lVert b\rVert_{\dot{H}^{\gamma - \alpha}}^{2}\lVert \nabla u\rVert_{\dot{H}^{\alpha}}^{2} \lesssim \lVert b\rVert_{\dot{H}^{\gamma - \alpha}}^{2}\lVert w\rVert_{\dot{H}^{\alpha}}^{2}.
\end{eqnarray*}
We then use Gagliardo-Nirenberg inequality, (7) and Proposition 3.1 to further bound by 

\begin{eqnarray*}
\lVert (b\cdot\nabla) u\rVert_{\dot{H}^{\gamma -1}}^{2} \lesssim (\lVert b\rVert_{L^{2}}^{2(1-(\gamma - \alpha))}\lVert \nabla b\rVert_{L^{2}}^{2(\gamma - \alpha)})\lVert w\rVert_{\dot{H}^{\alpha}}^{2} \lesssim \lVert w\rVert_{\dot{H}^{\alpha}}^{2}.
\end{eqnarray*}

Next, we fix $\epsilon \in (0, 1-\alpha)$ and estimate 

\begin{eqnarray*}
\lVert (u\cdot\nabla) b\rVert_{\dot{H}^{\gamma -1}}^{2}
&\lesssim &\lVert u\rVert_{\dot{H}^{\gamma - 1 + \epsilon}}^{2}\lVert \nabla b\rVert_{\dot{H}^{1-\epsilon}}^{2}\\
&\lesssim& \lVert u\rVert_{L^{2}}^{2(2-\gamma - \epsilon)}\lVert u\rVert_{\dot{H}^{1}}^{2(\gamma - 1 + \epsilon)}\lVert j\rVert_{\dot{H}^{1-\epsilon}}^{2} \lesssim \lVert j\rVert_{L^{2}}^{2\epsilon} \lVert j\rVert_{\dot{H}^{1}}^{2(1-\epsilon)} \lesssim (1+\lVert j\rVert_{\dot{H}^{1}}^{2})
\end{eqnarray*}
by Lemma 2.4, Gagliardo-Nirenberg inequalities, Lemma 2.1, Proposition 3.1 and Young's inequality. Thus, absorbing diffusive term, we have 

\begin{equation*}
\partial_{t} \lVert \Lambda^{\gamma} b\rVert_{L^{2}}^{2} +  \lVert \Lambda^{1+\gamma} b\rVert_{L^{2}}^{2} \lesssim (\lVert w\rVert_{\dot{H}^{\alpha}}^{2} + 1 + \lVert j\rVert_{\dot{H}^{1}}^{2}).
\end{equation*}
Hence, by Proposition 3.1, integrating in time we obtain 

\begin{equation*}
\sup_{t\in [0,T]} \lVert \Lambda^{\gamma} b(t)\rVert_{L^{2}}^{2} + \int_{0}^{T} \lVert \Lambda^{1+\gamma} b\rVert_{L^{2}}^{2} d\tau \leq c(u_{0}, b_{0}, T).
\end{equation*}
This completes the proof of Proposition 3.2.

}

\begin{proposition}

Let $N = 2, \nu, \eta > 0, \alpha \in (\frac{1}{3}, \sqrt{2} - 1), \beta = 1$. Then for any solution pair $(u,b)$ to (1) in $[0,T]$, for any $\gamma \in (1, 1+\alpha)$, there exists a constant $c(u_{0}, b_{0}, T)$ such that   

\begin{equation*}
\sup_{t\in [0,T]}\lVert w(t)\rVert_{L^{\frac{2(1+\alpha)}{2-\gamma}}}^{\frac{2(1+\alpha)}{2-\gamma}} + \int_{0}^{T} \lVert w\rVert_{L^{\frac{2(1+\alpha)}{(2-\gamma)(1-\alpha)}}}^{\frac{2(1+\alpha)}{2-\gamma}} d\tau \leq c(u_{0}, b_{0}, T).
\end{equation*}

\end{proposition}

\proof{

We fix $\gamma \in (1, 1+\alpha)$ and denote by 

\begin{equation*}
p = \frac{2(1+\alpha)}{2-\gamma}.
\end{equation*}
We estimate by multiplying the vorticity equation of (8) by $\lvert w\rvert^{p-2}w$ and integrating in space

\begin{eqnarray*}
\frac{1}{p} \partial_{t} \lVert w\rVert_{L^{p}}^{p} + \int \Lambda^{2\alpha} w\lvert w\rvert^{p-2} w dx = \int (b\cdot\nabla) j\lvert w\rvert^{p-2} w dx,
\end{eqnarray*}
where we used incompressibility of $u$. Using Lemma 2.3, because $p \geq 2$, and homogeneous Sobolev embedding $\dot{H}^{\alpha} \hookrightarrow L^{\frac{2}{1-\alpha}}$ we can obtain 

\begin{eqnarray*}
\int \Lambda^{2\alpha} w\lvert w\rvert^{p-2} w dx \geq \frac{2}{p} \lVert \lvert w\rvert^{\frac{p}{2}} \rVert_{\dot{H}^{\alpha}}^{2} \geq c(p,\alpha) \lVert w\rVert_{L^{\frac{p}{1-\alpha}}}^{p}.
\end{eqnarray*}
Using this, we further estimate 

\begin{eqnarray*}
\frac{1}{p}\partial_{t} \lVert w\rVert_{L^{p}}^{p} + c(p,\alpha) \lVert w\rVert_{L^{\frac{p}{1-\alpha}}}^{p} \leq  \lVert b\rVert_{L^{\infty}} \lVert \nabla j\rVert_{L^{\frac{2}{2-\gamma}}} \lVert w\rVert_{L^{p}}^{p-2} \lVert w\rVert_{L^{\frac{p}{1-\alpha}}},
\end{eqnarray*}
where we used the H$\ddot{o}$lder's inequality. Now we use the homogeneous Sobolev embedding of $\dot{H}^{\gamma - 1} \hookrightarrow L^{\frac{2}{2-\gamma}}$ and Gagliardo-Nirenberg inequality to obtain 

\begin{eqnarray*}
\frac{1}{p}\partial_{t} \lVert w\rVert_{L^{p}}^{p} + c(p,\alpha) \lVert w\rVert_{L^{\frac{p}{1-\alpha}}}^{p}
&\lesssim& \lVert b\rVert_{L^{\infty}} \lVert \Lambda^{\gamma} j\rVert_{L^{2}} \lVert w\rVert_{L^{p}}^{p-2} \lVert w\rVert_{L^{\frac{p}{1-\alpha}}}\\
&\lesssim& \lVert b\rVert_{L^{2}}^{\frac{\gamma - 1}{\gamma}}\lVert \Lambda^{\gamma} b\rVert_{L^{2}}^{\frac{1}{\gamma}} \lVert \Lambda^{\gamma} j\rVert_{L^{2}} \lVert w\rVert_{L^{p}}^{p-2} \lVert w\rVert_{L^{\frac{p}{1-\alpha}}}.
\end{eqnarray*}
We further bound by (7) and Proposition 3.2 to obtain 

\begin{eqnarray*}
&&\frac{1}{p}\partial_{t} \lVert w\rVert_{L^{p}}^{p} + c(p,\alpha) \lVert w\rVert_{L^{\frac{p}{1-\alpha}}}^{p}\\
&\leq& \frac{c(p,\alpha)}{2}\lVert w\rVert_{L^{\frac{p}{1-\alpha}}}^{p} + c\lVert \Lambda^{\gamma} j\rVert_{L^{2}}^{\frac{p}{p-1}}\lVert w\rVert_{L^{p}}^{(p-2)(\frac{p}{p-1})}\\
&\leq& \frac{c(p,\alpha)}{2} \lVert w\rVert_{L^{\frac{p}{1-\alpha}}}^{p} + c(1+\lVert \Lambda^{\gamma} j\rVert_{L^{2}}^{2})(1+ \lVert w\rVert_{L^{p}}^{p})
\end{eqnarray*}
by Young's inequalities. After absorbing the dissipative term, integrating in time gives 

\begin{equation*}
\sup_{t\in [0,T]} \lVert w(t)\rVert_{L^{p}}^{p} + \int_{0}^{T} \lVert w\rVert_{L^{\frac{p}{1-\alpha}}}^{p} d\tau \leq c(u_{0}, b_{0}, T)
\end{equation*}
due to Proposition 3.2. This completes the proof of Proposition 3.3. 

}

\begin{proposition}
Let $N = 2, \nu, \eta > 0, \alpha \in (\frac{1}{3}, \sqrt{2} - 1), \beta = 1$. Then for any solution pair $(u,b)$ to (1) in $[0,T]$, there exists a constant $c(u_{0}, b_{0}, T)$ such that   

\begin{equation*}
\sup_{t\in [0,T]}\lVert \nabla w(t)\rVert_{L^{2}}^{2} + \lVert \nabla j(t)\rVert_{L^{2}}^{2} + \int_{0}^{T}\lVert \Lambda^{\alpha} \nabla w\rVert_{L^{2}}^{2} + \lVert \Lambda \nabla j\rVert_{L^{2}}^{2} d\tau \leq c(u_{0}, b_{0}, T).
\end{equation*}

\end{proposition} 

\proof{

We apply $\nabla$ on (8) and take $L^{2}$-inner products with $\nabla w$ and $\nabla j$ respectively to estimate 

\begin{eqnarray*}
&&\frac{1}{2} \partial_{t} (\lVert \nabla w(t)\rVert_{L^{2}}^{2} + \lVert \nabla j(t)\rVert_{L^{2}}^{2}) +  \lVert \Lambda^{\alpha} \nabla w\rVert_{L^{2}}^{2} + \lVert \Lambda \nabla j\rVert_{L^{2}}^{2}\\
&=& -\int \nabla w \cdot \nabla u\cdot \nabla w dx - \int \nabla j \cdot \nabla u\cdot \nabla j dx\\
&&+ \int \nabla ((b\cdot\nabla)j)\cdot \nabla w + \nabla ((b\cdot\nabla)w)\cdot \nabla j dx\\
&&+ 2\int \nabla [\partial_{1}b_{1} (\partial_{1} u_{2} + \partial_{2} u_{1})]\cdot\nabla j dx - 2 \int \nabla [\partial_{1}u_{1}(\partial_{1}b_{2} + \partial_{2}b_{1})]\cdot \nabla j dx = \sum_{i=1}^{5}I_{i}.
\end{eqnarray*}

We estimate separately:

\begin{eqnarray*}
I_{1} &\leq& \lVert \nabla w\rVert_{L^{2}} \lVert \nabla u\rVert_{L^{\frac{2}{\alpha}}}\lVert \nabla w\rVert_{L^{\frac{2}{1-\alpha}}}\\
&\lesssim& \lVert \nabla w\rVert_{L^{2}} \lVert w\rVert_{L^{\frac{2}{\alpha}}}\lVert \Lambda^{\alpha} \nabla w\rVert_{L^{2}} \leq \frac{1}{8} \lVert \Lambda^{\alpha} \nabla w\rVert_{L^{2}}^{2} + cY(t) \lVert w\rVert_{L^{\frac{2}{\alpha}}}^{2}
\end{eqnarray*}
by H$\ddot{o}$lder's inequality, homogeneous Sobolev embedding of $\dot{H}^{\alpha} \hookrightarrow L^{\frac{2}{1-\alpha}}$, Lemma 2.1 and Young's inequality. Next, 

\begin{eqnarray*}
I_{2} \leq \lVert \nabla j\rVert_{L^{4}}^{2} \lVert \nabla u\rVert_{L^{2}}
\lesssim \lVert \nabla j\rVert_{L^{2}} \lVert \Delta j\rVert_{L^{2}} \lVert w\rVert_{L^{2}}
\leq \frac{1}{8} \lVert \Delta j\rVert_{L^{2}}^{2} + cY(t)
\end{eqnarray*}
by H$\ddot{o}$lder's, Gagliardo-Nirenberg and Young's inequalities. Next, we first integrate by parts and use the incompressibility conditions to obtain 

\begin{eqnarray*}
I_{3} &=& \int \nabla ((b\cdot\nabla)j) \cdot \nabla w dx + \int \nabla ((b\cdot\nabla)w)\cdot \nabla j dx\\
&=& \int \nabla b\cdot \nabla j \cdot \nabla w + \nabla b\cdot \nabla w \cdot \nabla j + b\cdot \nabla(\nabla j \otimes \nabla w) dx\\
&=& -\int \Delta b\cdot (\nabla j) w + 2\nabla b \cdot (\nabla \nabla j) w dx.
\end{eqnarray*}
We now estimate this by

\begin{eqnarray*}
I_{3} &\lesssim& (\lVert \nabla j\rVert_{L^{4}}^{2} \lVert w\rVert_{L^{2}} + \lVert \nabla b\rVert_{L^{4}} \lVert \Delta j\rVert_{L^{2}} \lVert w\rVert_{L^{4}})\\
&\lesssim& (\lVert \nabla j\rVert_{L^{2}} \lVert \Delta j\rVert_{L^{2}} \lVert w\rVert_{L^{2}} + \lVert j\rVert_{L^{2}}^{\frac{1}{2}} \lVert \nabla j\rVert_{L^{2}}^{\frac{1}{2}} \lVert \Delta j\rVert_{L^{2}} \lVert w\rVert_{L^{2}}^{\frac{1}{2}} \lVert \nabla w\rVert_{L^{2}}^{\frac{1}{2}})\\
&\leq& \frac{1}{8} \lVert \Delta j\rVert_{L^{2}}^{2} + cY(t)
\end{eqnarray*}
due to the H$\ddot{o}$lder's inequalities, Lemma 2.1, Gagliardo-Nirenberg inequalities, Proposition 3.1 and Young's inequalities. Finally, after integrating by parts again, 

\begin{eqnarray*}
I_{4} + I_{5} \lesssim \int \lvert \nabla b\rvert \lvert \nabla u\rvert \lvert \Delta j\rvert dx
\lesssim \lVert \nabla b\rVert_{L^{4}} \lVert w\rVert_{L^{4}} \lVert \Delta j\rVert_{L^{2}}
\end{eqnarray*}
by H$\ddot{o}$lder's inequality and Lemma 2.1. Note this is same as the second term of $I_{3}$ and hence its identical estimate suffices. 

Therefore, absorbing dissipative and diffusive terms, we have 

\begin{eqnarray*}
\partial_{t} Y(t) + \lVert \Lambda^{\alpha} \nabla w\rVert_{L^{2}}^{2} + \lVert \Lambda \nabla j\rVert_{L^{2}}^{2} \lesssim Y(t) (1+ \lVert w\rVert_{L^{\frac{2}{\alpha}}}^{2}).
\end{eqnarray*}

Now it can be checked that 

\begin{equation}
1 < 2-\frac{\alpha(1+\alpha)}{1-\alpha} < 1+\alpha \hspace{5mm} \forall \alpha \in (\frac{1}{3}, \sqrt{2} - 1)
\end{equation}
and hence we can choose $\gamma = 2-\frac{\alpha(1+\alpha)}{1-\alpha}$ so that by H$\ddot{o}$lder's inequality and Proposition 3.3,

\begin{equation*}
\int_{0}^{T} \lVert w\rVert_{L^{\frac{2}{\alpha}}}^{2} d\tau = \int_{0}^{T} \lVert w\rVert_{L^{\frac{2(1+\alpha)}{(2-\gamma)(1-\alpha)}}}^{2}d\tau \leq  T^{\frac{\gamma - 1 + \alpha}{1+\alpha}}\left(\int \lVert w\rVert_{L^{\frac{2(1+\alpha)}{(2-\gamma)(1-\alpha)}}}^{\frac{2(1+\alpha)}{2-\gamma}}\right)^{\frac{2-\gamma}{1+\alpha}} \leq c(u_{0}, b_{0}, T).
\end{equation*}
Therefore by Gronwall's inequality, 

\begin{equation*}
\sup_{t\in [0,T]} \lVert \nabla w(t)\rVert_{L^{2}}^{2} + \lVert \nabla j(t)\rVert_{L^{2}}^{2} + \int_{0}^{T} \lVert \Lambda^{\alpha} \nabla w\rVert_{L^{2}}^{2} + \lVert \Lambda \nabla j\rVert_{L^{2}}^{2} d\eta \leq c(u_{0}, b_{0}, T).
\end{equation*}
This completes the proof of Proposition 3.4. 

}

\textit{Proof of Theorem 1.1}

We now prove Theorem 1.1. We apply $\Lambda^{s}, s \in \mathbb{R}^{+}$ on (6) and take $L^{2}$-inner products with $\Lambda^{s}u$ and $\Lambda^{s}b$ respectively to estimate using Lemma 2.2 and incompressibility conditions to estimate

\begin{eqnarray*}
&&\partial_{t}(\lVert \Lambda^{s} u\rVert_{L^{2}}^{2} + \lVert \Lambda^{s} b\rVert_{L^{2}}^{2}) + \lVert \Lambda^{s + \alpha}u\rVert_{L^{2}}^{2} + \lVert \Lambda^{s + 1}b\rVert_{L^{2}}^{2}\\
&=& -\int \Lambda^{s} [(u\cdot\nabla)u]\cdot \Lambda^{s} u - u\cdot\nabla \Lambda^{s} u\cdot\Lambda^{s} u dx\\
&&- \int \Lambda^{s} [(u\cdot\nabla)b]\cdot\Lambda^{s} b - u\cdot\nabla \Lambda^{s} b\cdot\Lambda^{s} b dx\\
&&+ \int \Lambda^{s}[(b\cdot\nabla)b]\cdot\Lambda^{s} u - b\cdot\nabla \Lambda^{s} b\cdot\Lambda^{s} u dx\\
&&+ \int \Lambda^{s} [(b\cdot\nabla)u]\cdot\Lambda^{s} b - b\cdot\nabla\Lambda^{s} u\cdot\Lambda^{s} b dx\\
&\lesssim& (\lVert \nabla u\rVert_{L^{\frac{2}{\alpha}}} \lVert \Lambda^{s-1} \nabla u\rVert_{L^{2}} + \lVert \Lambda^{s} u\rVert_{L^{2}} \lVert \nabla u\rVert_{L^{\frac{2}{\alpha}}} ) \lVert \Lambda^{s} u\rVert_{L^{\frac{2}{1-\alpha}}}\\
&&+ (\lVert \nabla u\rVert_{L^{4}} \lVert \Lambda^{s-1} \nabla b\rVert_{L^{2}} + \lVert \Lambda^{s} u\rVert_{L^{2}} \lVert \nabla b\rVert_{L^{4}} ) \lVert \Lambda^{s} b\rVert_{L^{4}}\\
&&+ (\lVert \nabla b\rVert_{L^{4}} \lVert \Lambda^{s-1} \nabla b\rVert_{L^{4}} + \lVert \Lambda^{s} b\rVert_{L^{4}} \lVert \nabla b\rVert_{L^{4}}) \lVert \Lambda^{s} u\rVert_{L^{2}}\\
&&+ (\lVert \nabla b\rVert_{L^{4}} \lVert \Lambda^{s-1} \nabla u\rVert_{L^{2}} + \lVert \Lambda^{s} b\rVert_{L^{2}}\lVert \nabla u\rVert_{L^{4}}) \lVert \Lambda^{s} b\rVert_{L^{4}}\\
&\lesssim& (\lVert w\rVert_{L^{2}}^{\alpha}\lVert \nabla w\rVert_{L^{2}}^{1-\alpha}\lVert \Lambda^{s} u\rVert_{L^{2}})\lVert \Lambda^{s+\alpha} u\rVert_{L^{2}}\\
&&+(\lVert w\rVert_{L^{2}}^{\frac{1}{2}} \lVert \nabla w\rVert_{L^{2}}^{\frac{1}{2}} \lVert \Lambda^{s} b\rVert_{L^{2}} + \lVert \Lambda^{s} u\rVert_{L^{2}} \lVert j\rVert_{L^{2}}^{\frac{1}{2}} \lVert \nabla j\rVert_{L^{2}}^{\frac{1}{2}}) \lVert \Lambda^{s} b\rVert_{L^{2}}^{\frac{1}{2}} \lVert \Lambda^{s+1} b\rVert_{L^{2}}^{\frac{1}{2}}\\
&&+ (\lVert j\rVert_{L^{2}}^{\frac{1}{2}} \lVert \nabla j\rVert_{L^{2}}^{\frac{1}{2}} \lVert \Lambda^{s} b\rVert_{L^{2}}^{\frac{1}{2}} \lVert \Lambda^{s+1} b\rVert_{L^{2}}^{\frac{1}{2}}) \lVert \Lambda^{s} u\rVert_{L^{2}}\\
&&+ (\lVert j\rVert_{L^{2}}^{\frac{1}{2}} \lVert \nabla j\rVert_{L^{2}}^{\frac{1}{2}} \lVert \Lambda^{s} u\rVert_{L^{2}} + \lVert \Lambda^{s} b\rVert_{L^{2}} \lVert w\rVert_{L^{2}}^{\frac{1}{2}} \lVert \nabla w\rVert_{L^{2}}^{\frac{1}{2}})\lVert \Lambda^{s} b\rVert_{L^{2}}^{\frac{1}{2}} \lVert \Lambda^{s+1}b\rVert_{L^{2}}^{\frac{1}{2}}
\end{eqnarray*}
by H$\ddot{o}$lder's and Gagliardo-Nirenberg inequalities, homogeneous Sobolev embedding of $\dot{H}^{\alpha} \hookrightarrow L^{\frac{2}{1-\alpha}}$. Due to Propositions 3.1 and 3.4 and Young's inequalities we have 

\begin{eqnarray*}
&&\partial_{t}(\lVert \Lambda^{s} u\rVert_{L^{2}}^{2} + \lVert \Lambda^{s} b\rVert_{L^{2}}^{2}) + \lVert \Lambda^{s + \alpha}u\rVert_{L^{2}}^{2} + \lVert \Lambda^{s + 1}b\rVert_{L^{2}}^{2}\\
&\leq& \frac{1}{2} \left(\lVert \Lambda^{s+\alpha} u\rVert_{L^{2}}^{2} + \lVert \Lambda^{s+1} b\rVert_{L^{2}}^{2}\right) + c(\lVert \Lambda^{s} b\rVert_{L^{2}}^{2} +\lVert \Lambda^{s} u\rVert_{L^{2}}^{2}).
\end{eqnarray*}
Absorbing the dissipative and diffusive terms, Gronwall's inequality implies the desired result. 

\section{Proof of Theorem 1.2}

Throughout this section, we let $\alpha, \beta$ satisfy (4) and in particular we assume 

\begin{equation}
2\beta + \frac{\alpha(1+\alpha)}{1-\alpha} < 3 < \alpha + 2\beta + \frac{\alpha(1+\alpha)}{1-\alpha}
\end{equation}
as the other case can be done similarly. We work on 

\begin{equation}
\begin{cases}
\partial_{t}u + (u\cdot\nabla) u - (b\cdot\nabla) b + \nabla \pi + \Lambda^{2\alpha} u = 0,\\
\partial_{t}b + (u\cdot\nabla) b - (b\cdot\nabla) u + \Lambda^{2\beta} b = 0.
\end{cases}
\end{equation}

As before, taking $L^{2}$-inner products of (11) with $u$ and $b$ respectively, we immediately obtain 

\begin{equation}
\sup_{t\in [0,T]} \lVert u(t)\rVert_{L^{2}}^{2} + \lVert b(t)\rVert_{L^{2}}^{2} + \int_{0}^{T} \lVert \Lambda^{\alpha} u\rVert_{L^{2}}^{2} + \lVert \Lambda^{\beta} b\rVert_{L^{2}}^{2} d\tau \leq c(u_{0}, b_{0}, T).
\end{equation}

Since $\beta \geq 1$, it is clear from the proof of Proposition 3.1 that its slight modification applied to the following system 

\begin{equation}
\begin{cases}
\partial_{t} w + \Lambda^{2\alpha} w = -(u\cdot\nabla) w + (b\cdot\nabla) j\\
\partial_{t} j + \Lambda^{2\beta} j = -(u\cdot\nabla) j + (b\cdot\nabla) w + 2[\partial_{1}b_{1} (\partial_{1}u_{2} + \partial_{2}u_{1}) - \partial_{1}u_{1} (\partial_{1}b_{2} + \partial_{2}b_{1})]
\end{cases}
\end{equation}
leads to the following result: 

\begin{proposition} 
Let $N = 2, \nu, \eta > 0, \alpha \in (0, \frac{1}{3}], \beta \in (1,\frac{3}{2}]$ satisfy (10). Then for any solution pair $(u,b)$ to (1) in $[0,T]$, there exists a constant $c(u_{0}, b_{0}, T)$ such that 

\begin{equation*}
\sup_{t\in [0,T]} \lVert w(t)\rVert_{L^{2}}^{2} + \lVert j(t)\rVert_{L^{2}}^{2} + \int_{0}^{T} \lVert \Lambda^{\alpha} w\rVert_{L^{2}}^{2} + \lVert \Lambda^{\beta} j\rVert_{L^{2}}^{2} d\tau \leq c(u_{0}, b_{0}, T).
\end{equation*}

\end{proposition} 

Now we prove the following proposition: 

\begin{proposition}
Let $N = 2, \nu, \eta > 0, \alpha \in (0, \frac{1}{3}], \beta \in (1, \frac{3}{2}]$ satisfy (10). Then for any solution pair $(u, b)$ to (1) in $[0,T]$, for any $\gamma \in (\beta, \alpha + \beta)$, there exists a constant $c(u_{0}, b_{0}, T)$ such that 

\begin{equation*}
\sup_{t\in [0,T]} \lVert \Lambda^{\gamma} b(t)\rVert_{L^{2}}^{2} + \int_{0}^{T} \lVert \Lambda^{\beta + \gamma} b\rVert_{L^{2}}^{2} d\tau \leq c(u_{0}, b_{0}, T).
\end{equation*}

\end{proposition}

\proof{

We fix $\gamma \in (\beta, \alpha + \beta)$. From the magnetic field equation of (11), we estimate after multiplying by $\Lambda^{2\gamma}b$ and integrating in space

\begin{eqnarray*}
&&\frac{1}{2} \partial_{t} \lVert \Lambda^{\gamma} b\rVert_{L^{2}}^{2} + \lVert \Lambda^{\beta+ \gamma} b\rVert_{L^{2}}^{2}\\
 &\leq& (\lVert (u\cdot\nabla) b\rVert_{\dot{H}^{\gamma- \beta}}\lVert \Lambda^{\beta + \gamma} b\rVert_{L^{2}} + \lVert (b\cdot\nabla) u\rVert_{\dot{H}^{\gamma -\beta}} \lVert \Lambda^{\beta +\gamma} b\rVert_{L^{2}})\\
 &\leq& \frac{1}{2} \lVert \Lambda^{\beta + \gamma} b\rVert_{L^{2}}^{2} + c(\lVert (u\cdot\nabla) b\rVert_{\dot{H}^{\gamma - \beta}}^{2} + \lVert (b\cdot\nabla) u\rVert_{\dot{H}^{\gamma - \beta}}^{2})
 \end{eqnarray*} 
by H$\ddot{o}$lder's and Young's inequalities. 
Now by Lemma 2.4 and Gagliardo-Nirenberg inequality, Lemma 2.1 and Proposition 4.1, we estimate 

\begin{eqnarray*}
\lVert (b\cdot\nabla) u\rVert_{\dot{H}^{\gamma - \beta}}
\lesssim \lVert b\rVert_{\dot{H}^{\gamma - \beta + 1-\alpha}}^{2} \lVert \nabla u\rVert_{\dot{H}^{\alpha}}^{2} \lesssim \lVert b\rVert_{L^{2}}^{2(\alpha + \beta - \gamma)}\lVert \nabla b\rVert_{L^{2}}^{2(1-(\alpha + \beta - \gamma))}\lVert w\rVert_{\dot{H}^{\alpha}}^{2}
\lesssim \lVert w\rVert_{\dot{H}^{\alpha}}^{2}.
\end{eqnarray*}

Next, we fix $\epsilon \in (\beta - 1, \beta - \alpha)$ and  estimate using Lemma 2.4 and Gagliardo-Nirenberg inequalities, (12), Proposition 4.1 and Young's inequality as follows:

\begin{eqnarray*}
\lVert (u\cdot\nabla) b\rVert_{\dot{H}^{\gamma - \beta}}^{2}
&\lesssim& \lVert u\rVert_{\dot{H}^{\gamma + 1 - 2\beta + \epsilon}}^{2} \lVert \nabla b\rVert_{\dot{H}^{\beta - \epsilon}}^{2}\\
&\lesssim& \lVert u\rVert_{L^{2}}^{2(2\beta - \gamma - \epsilon)}\lVert \nabla u\rVert_{L^{2}}^{2(1-(2\beta - \gamma - \epsilon))}\lVert j\rVert_{\dot{H}^{\beta - \epsilon}}^{2}\\
&\lesssim& \lVert j\rVert_{L^{2}}^{2(\frac{\epsilon}{\beta})}\lVert j\rVert_{\dot{H}^{\beta}}^{2(1-\frac{\epsilon}{\beta})} \lesssim 1 + \lVert j\rVert_{\dot{H}^{\beta}}^{2}.
\end{eqnarray*}

Therefore, we have shown 

\begin{equation*}
\frac{1}{2}\partial_{t} \lVert \Lambda^{\gamma} b\rVert_{L^{2}}^{2} + \lVert \Lambda^{\beta + \gamma} b\rVert_{L^{2}}^{2} \leq \frac{1}{2} \lVert \Lambda^{\beta + \gamma} b\rVert_{L^{2}}^{2} + c(\lVert w\rVert_{\dot{H}^{\alpha}}^{2} + 1 + \lVert j\rVert_{\dot{H}^{\beta}}^{2}).
\end{equation*}
Integrating in time and using Proposition 4.1 complete the proof of Proposition 4.2.  

}  

\begin{proposition}

Let $N = 2, \nu, \eta > 0, \alpha \in (0, \frac{1}{3}], \beta \in (1, \frac{3}{2}]$ satisfy (10). Then for any solution pair $(u, b)$ to (1) in $[0,T]$, for any $\gamma \in (\beta, \alpha + \beta)$, there exists a constant $c(u_{0}, b_{0}, T)$ such that 

\begin{equation*}
\sup_{t\in [0,T]} \lVert w(t)\rVert_{L^{\frac{2(1+\alpha)}{3-\beta - \gamma}}}^{\frac{2(1+\alpha)}{3-\beta - \gamma}} + \int_{0}^{T} \lVert w\rVert_{L^{\frac{2(1+\alpha)}{(3-\beta - \gamma)(1-\alpha)}}}^{\frac{2(1+\alpha)}{3-\beta - \gamma}}d\tau \leq c(u_{0}, b_{0}, T).
\end{equation*}

\end{proposition}

\proof{

We fix $\gamma \in (\beta, \alpha + \beta)$ and denote 

\begin{equation*}
p = \frac{2(1+\alpha)}{3-\beta-\gamma}
\end{equation*}
Note due to (10), we have $3-\beta - \gamma > 0$. We estimate by multiplying the vorticity equation of (13) by $\lvert w\rvert^{p-2} w$ and integrating in space, using Lemma 2.3 and the same homogeneous Sobolev embedding of $\dot{H}^{\alpha} \hookrightarrow L^{\frac{2}{1-\alpha}}$ as before to obtain 

\begin{eqnarray*}
&&\frac{1}{p} \partial_{t} \lVert w\rVert_{L^{p}}^{p} + c(p,\alpha) \lVert w\rVert_{L^{\frac{p}{1-\alpha}}}^{p} \leq \lVert b\rVert_{L^{\infty}} \lVert \nabla j\rVert_{L^{\frac{p}{1+\alpha}}}\lVert w\rVert_{L^{p}}^{p-2} \lVert w\rVert_{L^{\frac{p}{1-\alpha}}}
\end{eqnarray*}
by H$\ddot{o}$lder's inequality. By our choice of $p$, we see that we may continue our estimate by 

\begin{eqnarray*}
&&\frac{1}{p} \partial_{t} \lVert w\rVert_{L^{p}}^{p} + c(p,\alpha) \lVert w\rVert_{L^{\frac{p}{1-\alpha}}}^{p}\\
&\leq& \lVert b\rVert_{L^{\infty}} \lVert \nabla j\rVert_{L^{\frac{2}{3-\beta - \gamma}}}\lVert w\rVert_{L^{p}}^{p-2} \lVert w\rVert_{L^{\frac{p}{1-\alpha}}}\\
&\lesssim& \lVert b\rVert_{L^{2}}^{\frac{\gamma - 1}{\gamma}}\lVert \Lambda^{\gamma} b\rVert_{L^{2}}^{\frac{1}{\gamma}} \lVert \Lambda^{\beta + \gamma} b\rVert_{L^{2}} \lVert w\rVert_{L^{p}}^{p-2} \lVert w\rVert_{L^{\frac{p}{1-\alpha}}}\\
&\leq& \frac{c(p,\alpha)}{2} \lVert w\rVert_{L^{\frac{p}{1-\alpha}}}^{p} + c(\lVert \Lambda^{\beta + \gamma} b\rVert_{L^{2}}^{2} + 1)(\lVert w\rVert_{L^{p}}^{p} + 1)
\end{eqnarray*}
where we used the Gagliardo-Nirenberg inequality, homogeneous Sobolev embedding of $\dot{H}^{\beta + \gamma - 2} \hookrightarrow L^{\frac{2}{3-\beta-\gamma}}$, Propositions 4.1 and 4.2, Young's inequalities. 

Absorbing dissipative term, Gronwall's inequality and Proposition 4.2 complete the proof of Proposition 4.3.

}

\begin{proposition}
Let $N  =2, \nu, \eta > 0, \alpha \in (0, \frac{1}{3}], \beta \in (1, \frac{3}{2}]$ satisfy (10). Then for any solution pair $(u,b)$ to (1) in $[0,T]$, there exists a constant $c(u_{0}, b_{0}, T)$ such that 

\begin{equation*}
\sup_{t\in [0,T]} \lVert \nabla w(t)\rVert_{L^{2}}^{2} + \lVert \nabla j(t)\rVert_{L^{2}}^{2} + \int_{0}^{T} \lVert \Lambda^{\alpha} \nabla w\rVert_{L^{2}}^{2} + \lVert \Lambda^{\beta} \nabla j\rVert_{L^{2}}^{2} d\tau \leq c(u_{0}, b_{0}, T).
\end{equation*}

\end{proposition}

\proof{

Similarly as before, we apply $\nabla$ on (13), take $L^{2}$-inner products with $\nabla w, \nabla j$ respectively to estimate 

\begin{eqnarray*}
&&\frac{1}{2} \partial_{t} (\lVert \nabla w(t)\rVert_{L^{2}}^{2} + \lVert \nabla j(t)\rVert_{L^{2}}^{2}) +  \lVert \Lambda^{\alpha} \nabla w\rVert_{L^{2}}^{2} + \lVert \Lambda^{\beta} \nabla j\rVert_{L^{2}}^{2}\\
&=& -\int \nabla w \cdot \nabla u\cdot \nabla w dx - \int \nabla j \cdot \nabla u\cdot \nabla j dx\\
&&+ \int \nabla ((b\cdot\nabla)j)\cdot \nabla w + \nabla ((b\cdot\nabla)w)\cdot \nabla j dx\\
&&+ 2\int \nabla [\partial_{1}b_{1} (\partial_{1} u_{2} + \partial_{2} u_{1})]\cdot\nabla j dx - 2 \int \nabla [\partial_{1}u_{1}(\partial_{1}b_{2} + \partial_{2}b_{1})]\cdot \nabla j dx = \sum_{i=1}^{5}I_{i}.
\end{eqnarray*}

As before,

\begin{eqnarray*}
I_{1} \leq \lVert \nabla w\rVert_{L^{2}} \lVert \nabla u\rVert_{L^{\frac{2}{\alpha}}} \lVert \nabla w\rVert_{L^{\frac{2}{1-\alpha}}}
\leq \frac{1}{8} \lVert \Lambda^{\alpha} \nabla w\rVert_{L^{2}}^{2} + cY(t) \lVert w\rVert_{L^{\frac{2}{\alpha}}}^{2}
\end{eqnarray*}
by H$\ddot{o}$lder's inequality, homogeneous Sobolev embedding of $\dot{H}^{\alpha} \hookrightarrow L^{\frac{2}{1-\alpha}}$ and Young's inequalities. Next, 

\begin{eqnarray*}
I_{2} 
\leq \lVert \nabla j\rVert_{L^{4}}^{2} \lVert \nabla u\rVert_{L^{2}}
\lesssim \lVert \nabla j\rVert_{L^{2}}^{2(\frac{2\beta - 1}{2\beta})}\lVert \Lambda^{\beta} \nabla j\rVert_{L^{2}}^{2(\frac{1}{2\beta})}
\leq \frac{1}{8} \lVert \Lambda^{\beta} \nabla j\rVert_{L^{2}}^{2}+ cY(t)
\end{eqnarray*}
by H$\ddot{o}$lder's inequality, Proposition 4.1, Gagliardo-Nirenberg and Young's inequalities. Next, we estimate $I_{3}$ after same integration by parts in the proof of Proposition 3.4,  

\begin{eqnarray*}
I_{3} &\lesssim& (\lVert \nabla j\rVert_{L^{4}}^{2} \lVert w\rVert_{L^{2}} + \lVert \nabla b\rVert_{L^{\frac{2}{\beta-1}}} \lVert \Delta j\rVert_{L^{\frac{2}{2-\beta}}}\lVert w\rVert_{L^{2}})\\
&\lesssim& \lVert \nabla j\rVert_{L^{2}}^{2(\frac{2\beta -1}{2\beta})}\lVert \Lambda^{\beta} \nabla j\rVert_{L^{2}}^{2(\frac{1}{2\beta})} + \lVert \nabla b\rVert_{L^{2}}^{\frac{2\beta + \gamma - 3}{\beta + \gamma - 1}}\lVert \Lambda^{\beta + \gamma}b\rVert_{L^{2}}^{\frac{2-\beta}{\beta + \gamma - 1}}\lVert \Lambda^{\beta} \nabla j\rVert_{L^{2}}\\
&\leq& \frac{1}{8} \lVert \Lambda^{\beta} \nabla j\rVert_{L^{2}}^{2} + c(Y(t) + 1 + \lVert \Lambda^{\beta + \gamma} b\rVert_{L^{2}}^{2})
\end{eqnarray*}
by H$\ddot{o}$lder's and Gagliardo-Nirenberg inequalities,  homogeneous Sobolev's embedding of $\dot{H}^{\beta - 1}\hookrightarrow L^{\frac{2}{2-\beta}}$ and Proposition 4.1. 

The estimates of $I_{4}$ and $I_{5}$ are simple: after the same integration by parts as before, we have 

\begin{eqnarray*}
I_{4} + I_{5} \lesssim \int \lvert \nabla b\rvert \lvert \nabla u\rvert \lvert \Delta j\rvert dx
\lesssim \lVert \nabla b\rVert_{L^{\frac{2}{\beta -1}}}\lVert \Delta j\rVert_{L^{\frac{2}{2-\beta}}} \lVert w\rVert_{L^{2}} 
\end{eqnarray*}
by H$\ddot{o}$lder's inequality and hence the same estimate as the second term of $I_{3}$ suffices. In sum, after absorbing dissipative and diffusive terms, we have  

\begin{eqnarray*}
&&\partial_{t} Y(t) + \lVert \Lambda^{\alpha} \nabla w\rVert_{L^{2}}^{2} + \lVert \Lambda^{\beta} \nabla j\rVert_{L^{2}}^{2} \lesssim (Y(t) + 1)(1+\lVert w\rVert_{L^{\frac{2}{\alpha}}}^{2} + \lVert \Lambda^{\beta + \gamma}b\rVert_{L^{2}}^{2}).
\end{eqnarray*}

Now we see that we may choose $\gamma = 3-\beta - \frac{\alpha(1+\alpha)}{1- \alpha}$ so that 

\begin{equation}
\beta < \gamma < \alpha + \beta
\end{equation}
due to (10) and therefore, by H$\ddot{o}$lder's inequality we have 

\begin{equation*}
\int_{0}^{T} \lVert w\rVert_{L^{\frac{2}{\alpha}}}^{2}d\tau \leq T^{\frac{\beta + \gamma - 2 + \alpha}{1+\alpha}}\left(\int_{0}^{T} \lVert w\rVert_{L^{\frac{2(1+\alpha)}{(3-\beta - \gamma)(1-\alpha)}}}^{\frac{2(1+\alpha)}{3-\beta-\gamma}}d\tau \right)^{\frac{3-\beta - \gamma}{1+\alpha}} \leq c(u_{0}, b_{0}, T)
\end{equation*}
due to Proposition 4.3. Thus, Gronwall's inequality and Proposition 4.2 complete the proof of Proposition 4.4. 

}

\textit{Proof of Theorem 1.2}

We are now ready to complete the proof of Theorem 1.2. Similarly as before we apply $\Lambda^{s}, s \in \mathbb{R}^{+}$ on (11) and take $L^{2}$-inner products with $\Lambda^{s}u$ and $\Lambda^{s}b$ respectively to estimate using Lemma 2.2

\begin{eqnarray*}
&&\partial_{t}(\lVert \Lambda^{s} u\rVert_{L^{2}}^{2} + \lVert \Lambda^{s} b\rVert_{L^{2}}^{2}) + \lVert \Lambda^{s + \alpha}u\rVert_{L^{2}}^{2} + \lVert \Lambda^{s + \beta}b\rVert_{L^{2}}^{2}\\
&\lesssim& \lVert \nabla u\rVert_{L^{\frac{2}{\alpha}}} \lVert \Lambda^{s} u\rVert_{L^{2}} \lVert \Lambda^{s} u\rVert_{L^{\frac{2}{1-\alpha}}}\\
&&+ (\lVert \nabla u\rVert_{L^{\frac{2}{1-\alpha}}}\lVert \Lambda^{s-1} \nabla b\rVert_{L^{\frac{2}{\alpha}}} + \lVert \Lambda^{s} u\rVert_{L^{\frac{2}{1-\alpha}}}\lVert \nabla b\rVert_{L^{\frac{2}{\alpha}}}) \lVert \Lambda^{s} b\rVert_{L^{2}}\\
&&+ \lVert \nabla b\rVert_{L^{\frac{2}{\alpha}}}\lVert \Lambda^{s} b\rVert_{L^{2}}\lVert \Lambda^{s} u\rVert_{L^{\frac{2}{1-\alpha}}}\\
&&+ (\lVert \nabla b\rVert_{L^{\frac{2}{\alpha}}}\lVert \Lambda^{s-1} \nabla u\rVert_{L^{\frac{2}{1-\alpha}}} + \lVert \Lambda^{s} b\rVert_{L^{\frac{2}{\alpha}}} \lVert \nabla u\rVert_{L^{\frac{2}{1-\alpha}}} ) \lVert \Lambda^{s} b\rVert_{L^{2}}\\
&\lesssim& \lVert w\rVert_{L^{2}}^{\alpha}\lVert \nabla w\rVert_{L^{2}}^{1-\alpha} \lVert \Lambda^{s} u\rVert_{L^{2}} \lVert \Lambda^{s+\alpha} u\rVert_{L^{2}}\\
&&+(\lVert \Lambda^{\alpha} w\rVert_{L^{2}} \lVert \Lambda^{s} b\rVert_{L^{2}}^{\frac{\alpha + \beta -1}{\beta}} \lVert \Lambda^{s+\beta} b\rVert_{L^{2}}^{\frac{1-\alpha}{\beta}}  + \lVert \Lambda^{s+\alpha} u\rVert_{L^{2}}\lVert \nabla b\rVert_{L^{2}}^{\alpha} \lVert \nabla \nabla b\rVert_{L^{2}}^{1-\alpha})\lVert \Lambda^{s} b\rVert_{L^{2}}\\
&&+\lVert \nabla b\rVert_{L^{2}}^{\alpha} \lVert \nabla \nabla b\rVert_{L^{2}}^{1-\alpha} \lVert \Lambda^{s} b\rVert_{L^{2}} \lVert \Lambda^{s+\alpha} u\rVert_{L^{2}}\\
&&+ (\lVert \nabla b\rVert_{L^{2}}^{\alpha} \lVert \nabla \nabla b\rVert_{L^{2}}^{1-\alpha} \lVert \Lambda^{s+\alpha} u\rVert_{L^{2}} + \lVert \Lambda^{s} b\rVert_{L^{2}}^{\frac{\alpha + \beta -1}{\beta}} \lVert \Lambda^{s+\beta} b\rVert_{L^{2}}^{\frac{1-\alpha}{\beta}} \lVert \Lambda^{\alpha} \nabla u\rVert_{L^{2}})\lVert\Lambda^{s} b\rVert_{L^{2}}\\
&\lesssim&  \lVert \Lambda^{s} u\rVert_{L^{2}} \lVert \Lambda^{s+\alpha} u\rVert_{L^{2}}\\
&&+ \lVert w\rVert_{L^{2}}^{1-\alpha} \lVert \nabla w\rVert_{L^{2}}^{\alpha} \lVert \Lambda^{s} b\rVert_{L^{2}}^{\frac{\alpha + 2\beta - 1}{\beta}} \lVert \Lambda^{s+\beta} b\rVert_{L^{2}}^{\frac{1-\alpha}{\beta}}  + \lVert \Lambda^{s+\alpha} u\rVert_{L^{2}} \lVert \Lambda^{s} b\rVert_{L^{2}}\\
&&+ \lVert \Lambda^{s} b\rVert_{L^{2}} \lVert \Lambda^{s+\alpha} u\rVert_{L^{2}}\\
&&+ \lVert \Lambda^{s} b\rVert_{L^{2}}^{\frac{\alpha + 2\beta - 1}{\beta}}\lVert \Lambda^{s+\beta} b\rVert_{L^{2}}^{\frac{1-\alpha}{\beta}} \lVert w\rVert_{L^{2}}^{1-\alpha} \lVert \nabla w\rVert_{L^{2}}^{\alpha}\\
&\leq& \frac{1}{2} \left(\lVert \Lambda^{s+\alpha} u\rVert_{L^{2}}^{2} +  \lVert \Lambda^{s+\beta} b\rVert_{L^{2}}^{2}\right) + c(\lVert \Lambda^{s} u\rVert_{L^{2}}^{2} + \lVert \Lambda^{s} b\rVert_{L^{2}}^{2})
\end{eqnarray*}
by H$\ddot{o}$lder's inequalities, Lemma 2.1, homogeneous Sobolev embedding of $\dot{H}^{\alpha}\hookrightarrow L^{\frac{2}{1-\alpha}}$ and Gagliardo-Nirenberg and Young's inequalities. Absorbing the dissipative and diffusive terms, Gronwall's inequality complete the proof of Theorem 1.2.


\begin{thebibliography}{100}  
 
\addtolength{\leftmargin}{0.2in} 
\setlength{\itemindent}{-0.2in} 

\bibitem[1]{1} C. Cao, and J. Wu, \textit{Two regularity criteria for the 3D MHD equations}, J. Differential Equations, \textbf{248}, 9 (2010), 2263-2274.

\bibitem[2]{2} C. Cao, and J. Wu, \textit{Global regularity for the 2D MHD equations with mixed partial dissipation and magnetic diffusion}, Adv. Math., \textbf{226}, 2 (2011), 1803-1822.

\bibitem[3]{3} C. Cao, J. Wu, and B. Yuan, \textit{The 2D incompressible magnetohydrodynamics equations with only magnetic diffusion}, arXiv:1306.3629 [math.AP], 10, Dec. 2013. 

\bibitem[4]{4} Q. Chen, C. Miao, and Z. Zhang, \textit{The Beale-Kato-Majda criterion for the 3D magneto-hydrodynamics equations}, Comm. Math. Phys., \textbf{275}, 3 (2007), 861-872.

\bibitem[5]{5} A. C$\acute{o}$rdoba, and D. C$\acute{o}$rdoba,  \textit{A maximum principle applied to quasi-geostrophic equations}, Comm. Math. Phys., \textbf{249}, 3 (2004), 511-528.

\bibitem[6]{6} Y. Du, Y. Liu, and Z. Yao, \textit{Remarks on the blow-up criteria for three-dimensional ideal magnetohydrodynamics equations}, J. Math. Phys., \textbf{50}, 023507 (2009).

\bibitem[7]{7} J. Fan, H. Gao, and G. Nakamura, \textit{Regularity criteria for the generalized magnetohydrodynamic equations and the quasi-geostrophic equations}, Taiwanese J. Math., \textbf{15}, 3 (2011), 1059-1073.

\bibitem[8]{8} A. Hasegawa, \textit{Self-organization processes in continuous media}, Adv. Phys., \textbf{34}, 1 (1985), 1-42.

\bibitem[9]{9} C. He, and Z. Xin, \textit{On the regularity of weak solutions to the magnetohydrodynamic equations}, J. Differential Equations, \textbf{213}, 2 (2005), 234-254. 

\bibitem[10]{10} X. Jia, and Y. Zhou, \textit{Regularity criteria for the 3D MHD equations involving partial components}, Nonlinear Anal. Real World Appl., \textbf{13}, 1 (2012), 410-418.

\bibitem[11]{11} Q. Jiu, and J. Zhao, \textit{Global regularity of 2D generalized MHD equations with magnetic diffusion}, arXiv:1309.5819 [math.AP], 10, Dec. 2013. 

\bibitem[12]{12} N. Ju, \textit{The maximum principle and the global attractor for the dissipative 2D quasi-geostrophic equations}, Comm. Math. Phys., \textbf{255}, 1 (2005), 161-181.

\bibitem[13]{13} T. Kato, and G. Ponce, \textit{Commutator estimates and the Euler and Navier-Stokes equations}, Comm. Pure Appl. Math., \textbf{41}, 7 (1988), 891-907.

\bibitem[14]{14} A. J. Majda, and A. L. Bertozzi, \textit{Vorticity and incompressible flow}, Cambridge University Press, 2001.

\bibitem[15]{15} R. May, \textit{Global well-posedness for a modified dissipative surface quasi-geostrophic equation in the critical Sobolev space $H^{1}$}, J. Differential Equations, \textbf{250}, 1 (2011), 320-339.

\bibitem[16]{16} C. Miao, and L. Xue, \textit{On the regularity of a class of generalized quasi-geostrophic equations}, J. Differential Equations, \textbf{251}, 10 (2011), 2789-2821.

\bibitem[17]{17} H. Miura, \textit{Dissipative quasi-geostrophic equation for large initial data in the critical Sobolev space}, Comm. Math. Phys., \textbf{267} (2006), 141-157.

\bibitem[18]{18} H. Politano, A. Pouquet, and P.-L. Sulem, \textit{Current and vorticity dynamics in three-dimensional magnetohydrodynamic turbulence}, Phys. Plasmas, \textbf{2} (1995), 2931-2939.

\bibitem[19]{19} M. Sermange, and R. Temam, \textit{Some mathematical questions related to the MHD equations}, Comm. Pure Appl. Math., \textbf{36} (1983), 635-664.

\bibitem[20]{20} T. Tao, \textit{Global regularity for a logarithmically supercritical hyperdissipative Navier-Stokes equation}, Anal. PDE, \textbf{2}, 3 (2009), 361-366.

\bibitem[21]{21} C. V. Tran, X. Yu, and L. A. K. Blackbourn, \textit{Two-dimensional magnetohydrodynamic turbulence in the limits of infinite and vanishing magnetic Prandtl number}, J. Fluid Mech., \textbf{725} (2013), 195-215. 

\bibitem[22]{22} C. V. Tran, X. Yu, and Z. Zhai, \textit{On global regularity of 2D generalized magnetohydrodynamics equations}, J. Differential Equations, \textbf{254}, 10 (2013), 4194-4216.

\bibitem[23]{23} C. V. Tran, X. Yu, and Z. Zhai, \textit{Note on solution regularity of generalized magnetohydrodynamic equations with partial dissipation}, Nonlinear Anal., \textbf{85} (2013), 43-51.

\bibitem[24]{24} J. Wu, \textit{Global regularity for a class of generalized magnetohydrodynamic equations}, J. Math. Fluid Mech., \textbf{13}, 2 (2011), 295-305. 

\bibitem[25]{25} J. Wu, \textit{Regularity criteria for the generalized MHD equations} Comm. Partial Differential Equations, \textbf{33}, 2 (2008), 285-306.

\bibitem[26]{26} J. Wu, \textit{The generalized MHD equations}, J. Differential Equations, \textbf{195} (2003), 284-312.

\bibitem[27]{27} K. Yamazaki, \textit{On the global regularity of generalized Leray-alpha type models}, Nonlinear Anal., \textbf{75}, 2 (2012), 503-515.

\bibitem[28]{28} K. Yamazaki, \textit{Remarks on the regularity criteria of generalized MHD and Navier-Stokes systems}, J. Math. Phys., \textbf{54}, 011502 (2013).

\bibitem[29]{29} K. Yamazaki, \textit{Global regularity of logarithmically supercritical MHD system with zero diffusivity}, Appl. Math. Lett., \textbf{29} (2014), 46-51.

\bibitem[30]{30} K. Yamazaki, \textit{A remark on the global well-posedness of a modified critical quasi-geostrophic equation}, ArXiV. arXiv:1006.0253 [math.AP], 10, Dec. 2013. 

\bibitem[31]{31} K. Yamazaki, \textit{Remarks on the global regularity of two-dimensional magnetohydrodynamics system with zero dissipation}, Nonlinear Anal., \textbf{94} (2014), 194-205.  

\bibitem[32]{32} B. Yuan and L. Bai, \textit{Remarks on global regularity of 2D generalized MHD equations}, arXiv:1306.2190 [math.AP], 10, Dec. 2013.

\bibitem[33]{33} Y. Zhou, \textit{Regularity criteria for the generalized viscous MHD equations}, Ann. I. H. Poincar$\acute{e}$ - AN \textbf{24}, 3 (2007), 491-505.

\bibitem[34]{34} Y. Zhou, \textit{Remarks on regularities for the 3D MHD equations}, Discrete Contin. Dyn. Syst., \textbf{12}, 5 (2005), 881-886.

\end{thebibliography}
\end{document}